\documentclass[11pt,twoside, spanish]{article}
\usepackage{psfrag,babel}
\usepackage{latexsym}
\usepackage{amsmath,amsthm}
\usepackage{amsfonts}
\usepackage{amssymb}
\usepackage{graphicx}
\usepackage{color}
\usepackage{epsfig}

\begin{document}
\title{Consistencia, no-trivialidad y redundancia en matem\'{a}tica.}
 \author{Eleonora Catsigeras$^1$}
\maketitle
\begin{center}
 
$^1$Instituto de Matem\'{a}tica y Estad\'istica \lq\lq Rafael Laguardia\rq\rq \\ Universidad de la Rep\'{u}blica,\\
Montevideo, Uruguay\\
{Corrreo electr\'{o}nico: eleonora@fing.edu.uy}

\end{center}

%%%%%%%%%%%%%%%%%%%%%%%%%%%%%%%%%%%%%%%%%%%%%%%%%%%%%%%%%%%%%%%%%%%%%%%%%%
%%%%%%%%%%%%%%%%%%%%%%%%%%%%%%%%%%%%%%%%%%%%%%%%%%%%%%%%%%%%%%%%%%%%%%%%%%
\begin{abstract}
Exploramos los criterios racionales formales e informales de consistencia, no-tri\-via\-lidad y  redundancia en la investigaci\'{o}n matem\'{a}tica actual. Desarrollamos la discusi\'{o}n paradigm\'{a}tica analizando las diferentes concepciones de esos criterios, desde las l\'{o}gico-formales hasta las informales (pero a\'{u}n racionales). Ilustramos la discusi\'{o}n con ejemplos concretos extra\'{\i}dos de la actividad de investigaci\'{o}n matem\'{a}tica, particularmente la publicada en los \'{u}ltimos cincuenta a\~{n}os en la teor\'{\i}a matem\'{a}tica de los sistemas din\'{a}micos deterministas.
\end{abstract}

\small

\noindent \em MSC 2010: \em 00A30, 03A05. 

\noindent \em Palabras y frases clave: \em  Racionalidad formal e informal, filosof\'{\i}a de la matem\'{a}tica, consistencia, no-trivialidad, profundidad matem\'{a}tica, redundancia matem\'{a}tica.

\vspace{.3cm}

\begin{center}
{\bf Abstract }
\end{center}

We explore the rational, formal and non-formal criteria of consistency, non-triviality and redundancy in the mathematical research now a days. We develop a paradigmatic discussion by analysing  the different conceptions of those criteria, from the logic-formal ones to the non formal ones (but still rational criteria).We illustrate the discussion with concrete examples obtained form the mathematical reseach, particularly from the published results that were published in the last 50 years in the mathematical theory of deterministic dynamical systems.

\vspace{.3cm}
\noindent \em MSC 2010: \em 00A30, 03A05.  

\noindent \em Key words and phrases: \em Formal and non-formal rationality, philosophy of mathematics, consistency, non-triviality, mathematical depth, mathematical redundancy.

\normalsize

\section{Introducci\'{o}n}

 Hilary Putnam \cite{15} defendi\'{o} que en toda justificaci\'{o}n racional hay impl\'{\i}cito un juicio de relevancia previo. Trata la relevancia en t\'{e}rminos de un modelo te\'{o}rico de valores en la ciencia. Bajo esa premisa, exploraremos tres de las muchas concepciones de las que surgen los juicios de relevancia de enunciados en la matem\'{a}tica actual, concepciones \'{e}stas tanto formales como informales.

Por un lado, la exploraci\'{o}n a lo largo de este trabajo es metaparadigm\'{a}tica, t\'{\i}picamente filos\'{o}fica, en la que discutiremos algunas de las concepciones racionales previas a los juicios de valor que la investigaci\'{o}n matem\'{a}tica actual asigna a sus enunciados. Por otro  lado, nuestra exploraci\'{o}n y discusi\'{o}n est\'{a}n basadas en la exposici\'{o}n de ejemplos de la pr\'{a}ctica cotidiana de la investigaci\'{o}n matem\'{a}tica, que muestra c\'{o}mo se interpretan en los hechos dos diferentes concepciones filos\'{o}ficas de racionalidad, diferentes pero complementarias, la l\'{o}gico-formal y la informal.

Adoptaremos un m\'{e}todo filos\'{o}fico por el cual, para sostener  un argumento filos\'{o}fico racional relativo al quehacer de las ciencias, es necesario (o por lo menos muy conveniente):

\em 1) nutrir el argumento con exposici\'{o}n de ejemplos;

2) no especular sobre c\'{o}mo una palabra o concepto deber\'{\i}a ser usado, sino observar c\'{o}mo se usa y aprender de ello. \em

(Wittgenstein \cite{23}, citado en Cook \cite[ p. 445]{4})

La metodolog\'{\i}a de esta investigaci\'{o}n, acorde a su prop\'{o}sito, es multidisciplinaria. Comprende dos disciplinas que aspiramos integrar en este trabajo: la matem\'{a}tica y la filosof\'{\i}a de la ciencia. Tomamos consciencia que los m\'{e}todos de exploraci\'{o}n y los lenguajes de comunicaci\'{o}n espec\'{\i}ficos a ambas disciplinas son en la actualidad  muy dis\'{\i}miles. Entre otras razones, estas disimilitudes son debidas al crecimiento y diversificaci\'{o}n de ambas disciplinas, y a la diferenciaci\'{o}n creciente de los prop\'{o}sitos de investigaci\'{o}n. Sin embargo, la metodolog\'{\i}a interdisciplinaria que adoptamos aqu\'{\i} intenta combinar en forma coherente los dos alfabetos de formas y lenguajes, uno tomado de la matem\'{a}tica y otro de la filosof\'{\i}a de la ciencia.

 \subsection {Convenci\'{o}n.} Acordamos la siguiente convenci\'{o}n, usual entre los investigadores matem\'{a}ticos. La palabra \lq\lq enunciado\rq\rq \ , lejos de acotar el objeto de estudio a su mera forma de expresi\'{o}n, se referir\'{a} al contenido matem\'{a}tico de un teorema, una definici\'{o}n o un axioma, en s\'{\i} mismo - sobre el cual discutiremos las caracter\'{\i}sticas meta-matem\'{a}ticas - bajo cualquier formulaci\'{o}n concreta con la que aparezca ese contenido. Se refiere entonces, por ejemplo, al teorema intr\'{\i}nsecamente, es decir, a las ideas matem\'{a}ticas en s\'{\i} mismas que el teorema involucra o relaciona, y a esta relaci\'{o}n, independientemente de su lenguaje y f\'{o}rmulas, de su notaci\'{o}n, expresi\'{o}n o formas de comunicaci\'{o}n concretas.

 \subsection {Prop\'{o}sito y tesis a defender. }
  
  La noci\'{o}n de \lq\lq relevancia\rq\rq \  que adoptaremos es una noci\'{o}n de consenso social. Se desprende del uso  habitual de esa palabra o concepto en la pr\'{a}ctica de evaluaci\'{o}n de resultados nuevos en la comunidad de investigadores matem\'{a}ticos. En particular, hemos considerado algunos de los juicios de valor de publicaciones matem\'{a}ticas, realizados por cient\'{\i}ficos matem\'{a}ticos independientes de los autores, y publicados en Zentralblatt f\"{u}r Mathematik y en Mathematical Reviews of the A.M.S.

Generalmente en un juicio de relevancia, la pr\'{a}ctica usual de la comunidad de matem\'{a}ticos, lejos de considerar \'{u}nicamente criterios l\'{o}gico-formales y exactos, adopta e integra (muchas veces con predominancia sobre la formalidad l\'{o}gica) otros criterios que, aunque informales, son a\'{u}n racionales. Justamente, la racionalidad de estos criterios informales, adem\'{a}s de los l\'{o}gico-formales, y la de su integraci\'{o}n con estos en un juicio de valor de enunciados matem\'{a}ticos, es la tesis central que defenderemos en este texto.

Nos enfocaremos en el an\'{a}lisis, la discusi\'{o}n y la ilustraci\'{o}n mediante ejemplos concretos, de tres conceptos claves y previos a los juicios de relevancia en la pr\'{a}ctica actual de la investigaci\'{o}n matem\'{a}tica (expondremos al final de esta introducci\'{o}n los motivos de tal enfoque). En efecto, a lo largo de las diferentes secciones de este trabajo discutiremos  nociones que se aplican en los juicios de valor de contenidos matem\'{a}ticos, en lo relativo a:

$\bullet$ consistencia,

$\bullet$ no-trivialidad, y

$\bullet$ no-redundancia.

Aunque cuando el matem\'{a}tico declara que un teorema es \lq\lq relevante\rq\rq \  o \lq\lq bueno\rq\rq \  o \lq\lq valioso\rq\rq,  no solo est\'{a}  afirmando que es consistente, no-trivial y no-redundante, nos centraremos en el an\'{a}lisis de esas tres nociones porque son ellas, desde  los puntos de vista l\'{o}gico-formal, semi-formal e informal, los principales y primeros (en orden cronol\'{o}gico) criterios de relevancia que se tienen en cuenta en la pr\'{a}ctica usual de evaluaci\'{o}n de enunciados matem\'{a}ticos nuevos. Muy frecuentemente son ellos tres, los \'{u}nicos pre-requisitos necesarios para un juicio de valor posterior, bastante m\'{a}s complejo, que integra una combinaci\'{o}n graduada y no cuantitativa de esos y varios otros criterios racionales de relevancia. As\'{\i}, los tres mencionados son condiciones necesarias previas para una justificaci\'{o}n de valor de los enunciados matem\'{a}ticos, aunque normalmente no sean suficientes.

\section { Criterios de consistencia l\'{o}gico-formal y experimental.}

Para definir consistencia definamos antes de concepto de contradicci\'{o}n. Dos enunciados matem\'{a}- ticos son contradictorios cuando de uno de ellos se puede deducir la negaci\'{o}n del otro. M\'{a}s restrictivamente, seg\'{u}n (Tarski, \cite[ p. 20]{21}) dos sentencias son contradictorias si una es equivalente a la negaci\'{o}n de la otra. Observamos que la definici\'{o}n de Tarski es restrictiva al requerir equivalencia entre sentencias. Sin embargo, lo \'{u}nico que se necesita para la inconsistencia de una pareja de proposiciones matem\'{a}ticas es que la una implique la negaci\'{o}n de la otra, y no necesariamente que esta otra implique la negaci\'{o}n de la primera.

La negaci\'{o}n de un teorema es l\'{o}gicamente diferente de la afirmaci\'{o}n contraria al teorema.  En efecto, la afirmaci\'{o}n contraria a un teorema de la forma \lq\lq si se cumple A entonces se cumple B\rq\rq \  es la afirmaci\'{o}n \lq\lq si no se cumple A entonces no se cumple B\rq\rq \ . La verdad del teorema es independiente de la verdad o falsedad de su afirmaci\'{o}n contraria. Por lo tanto, es consistente la pareja de afirmaciones compuestas por un teorema y su proposici\'{o}n contraria. En cambio, ese mismo teorema que enuncia \lq\lq si se cumple A entonces se cumple B\rq\rq \ , tiene como negaci\'{o}n \lq\lq existe un ejemplo que cumple A y no cumple B\rq\rq \  (llamado contraejemplo). El teorema es verdadero (si no lo fuera, no se llamar\'{\i}a teorema) si y solo si su negaci\'{o}n es falsa. Existe inconsistencia entre el teorema y su negaci\'{o}n, la existencia del contraejemplo.

Enunciados ya probados y futuros enunciados en una subdisciplina matem\'{a}tica (es decir, bajo un mismo sistema de definiciones y/o axiomas) son mutuamente consistentes cuando de ning\'{u}n subconjunto de ellos se puede deducir dos enunciados contradictorios.  As\'{\i}, un nuevo enunciado, o un nuevo supuesto en una subdisciplina matem\'{a}tica, es consistente con el estado de la subdisciplina, si y solo si del conjunto de todos los enunciados previos ya probados en ella no se puede deducir la negaci\'{o}n del nuevo enunciado o del nuevo supuesto.

Tarski (1946, \cite[p. 135]{21}) define la consistencia global de una disciplina deductiva, en vez de definir la consistencia de un nuevo enunciado en relaci\'{o}n al conjunto de enunciados previamente aceptados. En efecto seg\'{u}n Tarski, una teor\'{\i}a deductiva es llamada consistente o no contradictoria si ninguna pareja de afirmaciones de la teor\'{\i}a est\'{a} compuesta por una afirmaci\'{o}n y su negada.

Sin embargo, esta definici\'{o}n global no es exactamente la misma que la que se usa en la matem\'{a}tica contempor\'{a}nea. En efecto, la matem\'{a}tica de hoy en d\'{\i}a admite (como criterio racional informal) que cualquiera de sus  subdisciplinas est\'{a} compuesta solamente por la cantidad finita de enunciados, mutuamente consistentes, aceptados por los matem\'{a}ticos que investigan en esa subdisciplina, y por los que se pueden deducir de aquellos en forma trivial (definiremos trivialidad en la pr\'{o}xima secci\'{o}n). Es decir, no componen la subdisciplina todos los enunciados verdaderos posibles que se puedan deducir en forma no trivial en el futuro, a partir de los ya aceptados, ni los que se puedan crear en el futuro bajo nuevos presupuestos o definiciones. La matem\'{a}tica denomina preguntas abiertas a los enunciados no trivialmente verdaderos ni trivialmente falsos que relacionan conceptos definidos en la subdisciplina. Pero no emite un juicio de consistencia de las preguntas abiertas en relaci\'{o}n a los enunciados ya aceptados de la subdisciplina.

Estamos afirmando que al incorporar nuevos enunciados (demostrados y aceptados) la teor\'{\i}a no se conserva la misma, es decir, no se mantiene inmutable. La teor\'{\i}a o subdisciplina de la matem\'{a}tica, seg\'{u}n la conciben los investigadores matem\'{a}ticos, evoluciona con el tiempo, es din\'{a}mica, aunque cada enunciado ya incorporado no cambie. Se modifica el conjunto de enunciados ya demostrados y aceptados, por incorporaci\'{o}n de los nuevos, pero no por sustracci\'{o}n de los anteriores. Por lo tanto, cambia tambi\'{e}n el conjunto de preguntas abiertas.  Seg\'{u}n la definici\'{o}n global de consistencia de Tarski,  podr\'{\i}a interpretarse que para juzgar la consistencia de una teor\'{\i}a, se requerir\'{\i}a a priori el conocimiento de la totalidad de sus enunciados de manera est\'{a}tica, para saber si la disciplina es consistente o no lo es. Esta interpretaci\'{o}n requerir\'{\i}a del conocimiento a priori de todos los posibles enunciados verdaderos de la teor\'{\i}a. En cambio, la intepretaci\'{o}n parcial y din\'{a}mica de consistencia eval\'{u}a en forma continua el conjunto finito y cambiante de enunciados verdaderos (ya demostrados) de la disciplina en cada momento.

Observamos que la evoluci\'{o}n de la disciplina no solo se debe a que se incorporan nuevos enunciados que aparecen como resultados de la investigaci\'{o}n matem\'{a}tica. Se debe tambi\'{e}n a que para que estos nuevos enunciados sean consistentes con el conjunto de los enunciados anteriores, frecuentemente requieren la introducci\'{o}n de supuestos o definiciones adicionales.

Veamos un ejemplo de consistencia que requiri\'{o} la introducci\'{o}n de conceptos nuevos. Dentro de la teor\'{\i}a erg\'{o}dica diferenciable, una subdisciplina muy espec\'{\i}fica desarrollada durante los \'{u}ltimos 45 a\~{n}os, pero que presenta a\'{u}n numerosos problemas abiertos, es la llamada teor\'{\i}a de Pesin, y como parte de \'{e}sta, la teor\'{\i}a de medidas SRB (i.e. \lq\lq Sinai-Ruelle-Bowen\rq\rq), originada entre otros, en el art\'{\i}culo precursor de Sinai en 1972, \cite{18}. Esas dos sub-ramas de la teor\'{\i}a erg\'{o}dica diferenciable asumen que el sistema din\'{a}mico bajo estudio tiene regularidad $C^2$: el sistema es diferenciable hasta segundo orden con derivadas hasta segundo orden continuas. La regularidad $C^2$ es necesaria  para las demostraciones conocidas de algunos teoremas muy populares en esas subdisciplinas. En efecto, en (Robinson-Young, \cite[pp. 159-176]{16}) se muestra un contraejemplo sin regularidad $C^2$, para el cual algunos resultados de la teor\'{\i}a de Pesin no rigen. El contraejemplo de Robinson y Young es solo de clase $C^1$: el sistema es diferenciable solamente hasta primer orden con derivada primera continua, pero no lo es hasta segundo orden. En particular, la existencia de medidas SRB\footnote{Estamos considerando las medidas SRB como las que satisfacen ciertas propiedades de continuidad absoluta. M\'{a}s precisamente,  la proyecci\'{o}n a lo largo de la foliaci\'{o}n estable y las medidas condicionales respecto a la foliaci\'{o}n inestable deben ser absolutamente continuas.} es falsa en ese contraejemplo.

Como contraparte, en (Enrich-Catsigeras, \cite[p. 740]{5}) se demuestra que todos los sistemas continuos, incluyendo los de clase $C^1$ y los de clase $C^2$, poseen cierto tipo de medidas abstractas llamadas \lq\lq SRB-like\rq\rq \  (i.e. parecidas a las medidas llamadas SRB), definidas como las que presentan la propiedad de observabilidad a trav\'{e}s de los promedios temporales asint\'{o}ticos de conjuntos Lebesgue-positivos de \'{o}rbitas del sistema. Esta propiedad de observabilidad es exhibida en particular por todas las medidas SRB de los sistemas diferenciables, en los casos en que las SRB existen. Es decir, todas las medidas SRB de los sistemas din\'{a}micos diferenciables son casos particulares de medidas SRB-like, y estas \'{u}ltimas siempre existen, aunque el sistema no sea diferenciable.

?`No es el nuevo enunciado de existencia de medidas SRB-like inconsistente con los contraejemplos conocidos previamente para los que no existen medidas SRB? La respuesta es \lq\lq no, no lo es\rq\rq, pues el nuevo teorema se refiere a las medidas SRB-like y no a las SRB solamente. \'{E}stas son, cuando existen, solo caso particulares de aquellas. Por lo tanto pueden no existir las SRB, como en el contraejemplo previamente conocido de  Robinson y Young, y s\'{\i} existir las SRB-like  como se demuestra en el nuevo teorema.

En conclusi\'{o}n, la consistencia del nuevo teorema con el cuerpo de la subdisciplina desarrollada hasta ese momento, requiri\'{o} de la creaci\'{o}n e introducci\'{o}n de un concepto matem\'{a}tico nuevo: el de las medidas SRB-like o propiedad de observabilidad. Este concepto nuevo es consistente con el concepto previo de medida SRB, y m\'{a}s a\'{u}n, aparentemente fue extra\'{\i}do o inspirado en el descubrimiento de la propiedad necesaria de observabilidad de \'{e}stas.

Discutamos ahora la concepci\'{o}n filos\'{o}fica de consistencia experimental, adem\'{a}s de la consistencia l\'{o}gico-formal que ya expusimos. Esta consideraci\'{o}n de conceptos filos\'{o}ficos diferentes de consistencia matem\'{a}tica (pero no opuestos sino complementarios), permite diversificar las estrategias para excluir el error. En matem\'{a}tica, la consistencia experimental es un criterio de relevancia de un nuevo teorema: significa que el teorema nuevo que se demuestra, no solo es consistente desde el punto de vista l\'{o}gico-formal con los enunciados aceptados hasta el momento en la subdisciplina matem\'{a}tica a la cual se incorpora, sino que explica (o demuestra la necesidad de) propiedades matem\'{a}ticas ya observadas previamente en el comportamiento de ejemplos o casos particulares no-triviales, que por alguna raz\'{o}n son considerados relevantes o paradigm\'{a}ticos dentro de la disciplina.

Es m\'{a}s f\'{a}cil encontrar ejemplos donde la consistencia experimental es relativa a la investigaci\'{o}n de la matem\'{a}tica aplicada a otras ciencias, que de la matem\'{a}tica pura.  As\'{\i}, por ejemplo, el teorema de (Mirollo-Strogatz, \cite[ pp. 1645-1662]{13}) en la subdisciplina matem\'{a}tica de los sistemas din\'{a}micos deterministas, (como rama de la bio-matem\'{a}tica, y en particular aplicada a la neurociencias),  demuestra la necesaria sincronizaci\'{o}n del evento llamado \lq\lq espiga\rq\rq \  o \lq\lq disparo\rq\rq \  de las neuronas de una red excitatoria. Para poder enunciar y demostrar ese teorema, est\'{a}n definidas las neuronas, la red que conforman, y el fen\'{o}meno de espiga, de forma matem\'{a}tica abstracta modelando la red o subred neural\footnote{Cuando la red est\'{a} definida en forma abstracta matem\'{a}tica, o dise\~{n}ada por la ingenier\'{\i}a como red artificial, se llama "red neuronal". Cuando es una red de neuronas biol\'{o}gicas del sistema nervioso de un animal se llama "red neural". } biol\'{o}gica concreta. Ese teorema, desde el punto de vista l\'{o}gico-formal de su enunciado y de su demostraci\'{o}n, es relevante porque entre otros motivos, es consistente con el fen\'{o}meno de sincronizaci\'{o}n de espigas observado experimentalmente mediante encefalogramas en ciertas subredes neurales biol\'{o}gicas excitatorias estudiadas previamente por los neurocient\'{\i}ficos. Por lo tanto, la relevancia por consistencia experimental del teorema de Mirollo y Strogatz, est\'{a} en parte fundado en la relevancia de los resultados previos obtenidos experimentalmente en esas otras disciplinas o ciencias, a las cuales el teorema se aplica.  Este criterio de relevancia de enunciados matem\'{a}ticos aplicados o aplicables, basado en la relevancia de resultados experimentales previos o posteriores, con los cuales el enunciado es consistente, no es un criterio que se limita a verificaciones l\'{o}gico-formales. Es racional informal, apelando al siguiente argumento  de Putnam sobre la racionalidad de las ciencias:

\em Los procedimientos mediante los que decidimos (la relevancia) tienen que ver con que considerada como un todo, exhiba ciertas 'virtudes', o no las exhiba. Estoy suponiendo que el procedimiento ... (de decisi\'{o}n), no puede analizarse correctamente como un procedimiento de verificaci\'{o}n ... oraci\'{o}n por oraci\'{o}n. Estoy suponiendo que la verificaci\'{o}n... es una cuesti\'{o}n hol\'{\i}stica, que son los sistemas... enteros los que se enfrentan al tribunal de la experiencia..., y que el juicio resultante es una cuesti\'{o}n un tanto intuitiva, que no puede formalizarse a menos que formalicemos toda la psicolog\'{\i}a humana... \rm
 (Putnam 1981, \cite[cap. 6, par. 1]{15})

Aunque es m\'{a}s dif\'{\i}cil encontrar ejemplos, en la investigaci\'{o}n matem\'{a}tica pura la consistencia experimental tambi\'{e}n es un criterio racional informal complementario al l\'{o}gico-formal, y aplicable a juicios de relevancia de enunciados nuevos. Por ejemplo, Anosov en 1967 \cite{1} defini\'{o} una clase de difeomorfismos, y dio ejemplos particulares de ellos, que son uniformemente hiperb\'{o}licos en forma global en todo el espacio. Estos difeomorfismos fueron posteriormente llamados \lq\lq difeomorfismos de Anosov\rq\rq \ .  Todos los ejemplos de esta clase de difeomorfismos estudiados hasta el d\'{\i}a de hoy (descontados aquellos ejemplos para los que a\'{u}n no se conocen las propiedades erg\'{o}dicas), cumplen la siguiente propiedad:

\begin{center}
\em 
Si un difeomorfismo de Anosov es conservativo\footnote{Un difeomorfismo se llama conservativo si preserva la forma de volumen en el espacio donde act\'{u}a.}, entonces es erg\'{o}dico\footnote{Un difeormofismo se llama erg\'{o}dico si los promedios temporales asint\'{o}ticos de funciones observables a lo largo de casi todas sus \'{o}rbitas, son iguales al promedio espacial o valor esperado de la funci\'{o}n que se observa.  }. 
\end{center}

\rm

El enunciado anterior no es un teorema. Es solamente el resultado obtenido en todos los ejemplos de difeomorfismos de Anosov conservativos para los cuales se pudo estudiar la ergodicidad. Es a\'{u}n una pregunta abierta si la proposici\'{o}n anterior es verdadera o falsa para la totalidad de los difeomorfismos de Anosov conservativos. Es decir, la ergodicidad de los difeomorfismos de Anosov  conservativos es un resultado \lq\lq experimental\rq\rq \  de la matem\'{a}tica pura; y lo llamamos \lq\lq experimental\rq\rq \  porque su verificaci\'{o}n, aunque sea l\'{o}gico-formal y abstracta, es aplicable solamente uno a uno, a casos particulares que no son la generalidad.

En el a\~{n}o 1972, Sinai \cite{18} demostr\'{o} que todos los difeomorfismos de Anosov conservativos  que  son de clase $C^2$,  son erg\'{o}dicos. Este teorema de Sinai, adem\'{a}s de ser relevante seg\'{u}n los criterios de no-trivialidad que discutiremos en la secci\'{o}n siguiente, lo es por su propiedad de consistencia experimental con todos los ejemplos de difeomorfismos de Anosov cuya ergodicidad era conocida en forma particular, caso a caso,  hasta el momento de su aparici\'{o}n. Y aunque no demuestra la ergodicidad en todos los ejemplos posibles (pues excluye los difeomorfismos que no son de clase $C^2$), extiende la propiedad de ergodicidad, antes conocida ejemplo a ejemplo solamente, a toda una subclase de difeomorfismos de Anosov conservativos.

Concluimos que el criterio de consistencia experimental de un nuevo teorema es un concepto racional  informal, ya que explica en general las propiedades exhibidas en ejemplos previamente considerados relevantes en esa subdisciplina; y la relevancia no es una propiedad verificable mediante un algoritmo de deducci\'{o}n, ni apelando a argumentos l\'{o}gico-formales exclusivamente.

En las pr\'{o}ximas secciones extenderemos las estrategias que hemos utilizado hasta ahora, que incluyen dos concepciones filos\'{o}ficas diferentes pero complementarias de  consistencia matem\'{a}tica (la l\'{o}gico-formal y la experimental),  para discutir tambi\'{e}n la no-trivialidad y la redundancia de los enunciados matem\'{a}ticos.

\section{Criterios de no-trivialidad l\'{o}gico-formal e informal.}

El adjetivo \lq\lq no-trivial\rq\rq \  es com\'{u}nmente usado por los matem\'{a}ticos para indicar que un teorema no es obvio o no es f\'{a}cil de demostrar (Ver por ejemplo Weisstein, \cite[p. 1]{22}). Aqu\'{\i}, es importante no sustituir  la conjunci\'{o}n \lq\lq o\rq\rq \  por la conjunci\'{o}n \lq\lq y\rq\rq \ . Es decir, un teorema puede ser no-trivial, si no es obvio lo que enuncia, aunque su demostraci\'{o}n sea breve y f\'{a}cil de entender.  Pero esta definici\'{o}n, requiere definir primero qu\'{e} significa obviedad, y qu\'{e} significa f\'{a}cil de demostrar.  Y estas propiedades de un teorema no son propiedades l\'{o}gico-formales, sino informales, aunque plausibles  de producir un criterio racional de relevancia del teorema.  Como dependen del contexto y del investigador (matem\'{a}tico o no), el definido m\'{a}s arriba es un criterio informal. Analizaremos m\'{a}s adelante este criterio. Es el que adoptamos como criterio de no-trivialidad de un teorema, justificando nuestra elecci\'{o}n a lo largo de la siguiente discusi\'{o}n.

Veamos varios posibles significados, estrechamente vinculados entre s\'{\i}, de la no-trivialidad l\'{o}gico-formal. Para ello definamos su negaci\'{o}n, la trivialidad. Una nueva afirmaci\'{o}n es trivial, desde el punto de vista estricto l\'{o}gico-formal, si de los enunciados (definiciones, axiomas, teoremas) conocidos en la subdisciplina matem\'{a}tica desarrollada hasta el momento, se puede deducir esa afirmaci\'{o}n. Esta definici\'{o}n no es nunca empleada por los matem\'{a}ticos investigadores. Si la aplicaran, todos sus resultados ser\'{\i}an triviales, independientemente de la profundidad de los significados, de las construcciones de nuevos conceptos matem\'{a}ticos involucradas en su obtenci\'{o}n, y de la dificultad y longitud de las demostraciones. Son estos factores, y no la definici\'{o}n estricta l\'{o}gico-formal, los decisivos en un juicio de trivialidad o no-trivialidad de enunciados para los matem\'{a}ticos.

Agreguemos entonces la siguiente condici\'{o}n: una nueva afirmaci\'{o}n es trivial, desde un punto de vista cuantitativo l\'{o}gico-formal,  si se puede deducir  mediante una cantidad breve de pasos deductivos directos. Para lo cual, es necesario definir a priori  la cota superior num\'{e}rica de la \lq\lq cantidad breve\rq\rq \  y la lista finita de los pasos deductivos \lq\lq elementales o directos\rq\rq \ .

En este sentido, una afirmaci\'{o}n es trivial para un matem\'{a}tico si su demostraci\'{o}n es un \lq\lq ejercicio\rq\rq \ . Efectivamente, el matem\'{a}tico puro que no trabaja en aplicaciones,   no publica como resultados de investigaci\'{o}n lo que considera ejercicios, pues los eval\'{u}a como no trascendentes a la investigaci\'{o}n matem\'{a}tica en su especialidad, y por lo tanto irrelevantes para el avance de \'{e}sta.  Pero, en sentido contrario, desde el punto de vista del matem\'{a}tico aplicado, y de los investigadores de otras ciencias que aplican la matem\'{a}tica, la postura de considerar triviales los enunciados brevemente deducibles a partir de lo ya conocido en la subdisciplina matem\'{a}tica que se aplica, resulta una limitaci\'{o}n inadecuada. En efecto, seg\'{u}n el premiado Nobel en f\'{\i}sica, Richard Feynman \cite{6},\lq\lq los matem\'{a}ticos designan cualquier teorema como trivial, una vez que su prueba ya fue obtenida y es conocida, ... Por lo tanto hay solo dos tipos de proposiciones matem\'{a}ticas verdaderas: las triviales y las a\'{u}n no demostradas\rq\rq \  (Feynman 1997, citado en Weisstein, \cite[p. 1]{22}).

Los criterios de \lq\lq no-trivialidad\rq\rq \  matem\'{a}tica que la mayor parte de los matem\'{a}ticos investigadores utilizan, ya sea en matem\'{a}tica pura o aplicada, no son \'{u}nicamente los l\'{o}gico-formales. Esta es la primera de las razones por las que preferimos afiliarnos al criterio racional informal de no-trivialidad, definido al principio de esta secci\'{o}n, y que discutiremos a continuaci\'{o}n.

Entre otros autores, Shanks \cite{17} define \lq\lq el opuesto de un teorema trivial ...(como) un 'teorema profundo'. Cualitativamente, un teorema profundo es un teorema cuya prueba es larga, complicada o dif\'{\i}cil, o (cuyo enunciado) involucra temas de la matem\'{a}tica que no est\'{a}n obviamente relacionadas...\rq\rq \  (Shanks, \cite[pp. 64-66]{17}) Esta definici\'{o}n de no-trivialidad matem\'{a}tica, como sin\'{o}nimo de profundidad del enunciado, en cuanto a la relaci\'{o}n entre conceptos matem\'{a}ticos no primariamente vinculados, es recogida y preferida en la actualidad por buena parte de los investigadores en matem\'{a}tica. Por ejemplo (Tao, \cite[p. 623-624]{20}) enumera veintid\'{o}s criterios no exhaustivos  para explicar qu\'{e} es la \lq\lq buena\rq\rq \  matem\'{a}tica. Entre esos criterios, el concepto de profundidad de un resultado matem\'{a}tico es considerado por Tao como una caracter\'{\i}stica de evidente no-trivialidad, porque, entre otras razones, capta un fen\'{o}meno m\'{a}s all\'{a} del alcance de herramientas m\'{a}s elementales.

Seg\'{u}n Gray \cite{8} \lq\lq los matem\'{a}ticos usan la palabra 'profundo' para referir una alta apreciaci\'{o}n de un concepto, teorema o demostraci\'{o}n.  Su primer uso en matem\'{a}ticas fue una consecuencia del trabajo de Gauss en teor\'{\i}a de n\'{u}meros y el acuerdo entre sus sucesores que partes espec\'{\i}ficas de ese trabajo mostraban propiedades estructurales de la matem\'{a}tica... en contraste al alcance menos estructural y m\'{a}s orientado a la soluci\'{o}n de problemas\rq\rq \  de los trabajos matem\'{a}ticos menos profundos  (Gray, \cite[p. 177]{8}).

La profundidad de un teorema est\'{a} referida al criterio racional informal de no-trivialidad, y a la vinculaci\'{o}n a trav\'{e}s de su enunciado de estructuras matem\'{a}ticas no directamente relacionados o evidentes en forma previa. Cuando esta vinculaci\'{o}n resulta  sorprendente o inesperada, a\'{u}n si la demostraci\'{o}n formal es breve y sencilla, el teorema es profundo.  Es un criterio informal, pues no solo es irreducible a algoritmos formales, sino que involucra al sujeto en su percepci\'{o}n de profundidad.  Aunque informal, es un criterio racional de no-trivialidad, por ejemplo cuando se la justifica y fundamenta matem\'{a}ticamente, de manera independiente de la demostraci\'{o}n formal del teorema. As\'{\i}, este criterio racional de no-trivialidad es inherente al enunciado mismo, en forma desligada de su demostraci\'{o}n. Puede ser adem\'{a}s, inmutable en el tiempo, aunque en el futuro puedan encontrarse demostraciones formales muy simples del mismo teorema.

Ilustremos el argumento anterior con un ejemplo que no corresponde a la investigaci\'{o}n de la matem\'{a}tica de hoy en d\'{\i}a, sino a la matem\'{a}tica de la Antig\"{u}edad. El teorema de Pit\'{a}goras es no-trivial, intr\'{\i}nsecamente, pues en forma inesperada y sorpresiva relaciona las operaciones de suma y multiplicaci\'{o}n de cantidades num\'{e}ricas con  las longitudes de los lados de un tri\'{a}ngulo rect\'{a}ngulo. Se han hecho populares hoy en d\'{\i}a demostraciones brev\'{\i}simas y muy sencillas del Teorema de Pit\'{a}goras.  Pero a\'{u}n as\'{\i} ?`c\'{o}mo pueden las operaciones num\'{e}ricas de multiplicaci\'{o}n y suma ser tan una manifestaci\'{o}n de la geometr\'{\i}a m\'{e}trica de un tri\'{a}ngulo rect\'{a}ngulo? Los conceptos est\'{a}n definidos, a priori, de forma  independiente. ?`Por qu\'{e} sucede entonces que esa aparente independencia no sea tal, sino que sean conceptos vinculados a trav\'{e}s del Teorema de Pit\'{a}goras? Entender las causas profundas de relaciones matem\'{a}ticas conlleva una comprensi\'{o}n diferente de la mera explicaci\'{o}n l\'{o}gica, deductiva y formal de cada paso de una demostraci\'{o}n, e inspira la creaci\'{o}n de definiciones matem\'{a}ticas abstractas que reflejan esa comprensi\'{o}n profunda. En (Livio, \cite[pp. 80-83]{11}) se sugiere que esa comprensi\'{o}n trasciende la forma l\'{o}gico-deductiva de relaciones no obvias, y en ella radica  la respuesta \lq\lq es ambos\rq\rq, a la vieja pregunta de si la matem\'{a}tica es invento o descubrimiento.

El criterio de no-trivialidad por no-obviedad o profundidad del enunciado de un teorema, a pesar de la eventual simplicidad de su demostraci\'{o}n formal, est\'{a} estrechamente ligado a los criterios racionales informales de relevancia de un enunciado matem\'{a}tico por su belleza est\'{e}tica y por su simplicidad. La mayor\'{\i}a de los  investigadores en la matem\'{a}tica pura reconocen su motivaci\'{o}n en la b\'{u}squeda de la no-trivialidad informal, es decir la  profundidad del significado de los enunciados que crean o descubren. Adem\'{a}s encuentran belleza y relevancia en el contraste entre la no-obviedad y la simplicidad de la relaci\'{o}n hallada, cuando la demostraci\'{o}n que encuentran, desde el punto de vista l\'{o}gico-formal, se puede reducir a una cantidad breve de pasos deductivos relativamente elementales.  Por ejemplo Hardy \cite{9} sugiere que una de las condiciones necesarias para la percepci\'{o}n de belleza en la matem\'{a}tica es la econom\'{\i}a de recursos formales.

Sin embargo, en la pr\'{a}ctica de la comunicaci\'{o}n y publicaci\'{o}n de los resultados de la matem\'{a}tica pura, resulta hoy en d\'{\i}a muy dif\'{\i}cil encontrar ejemplos de enunciados nuevos no-triviales, es decir profundos o no-obvios, pero con demostraciones breves y sencillas.  En primer lugar, el propio matem\'{a}tico frecuentemente auto-censura sus resultados como triviales, y solo los difunde en car\'{a}cter de ejercicios, cuando obtiene demostraciones breves y sencillas, aunque sus enunciados sean nuevos, y muestren propiedades inesperadas y sorprendentes. En segundo lugar, si el investigador matem\'{a}tico a\'{u}n cree relevante su resultado a pesar de la sencillez de su demostraci\'{o}n, y  lo somete a publicaci\'{o}n, dif\'{\i}cilmente el \'{a}rbitro se convencer\'{a} racionalmente de la no-trivialidad. Frecuentemente, en un juicio de no trivialidad se considera primero el aspecto formal de la demostraci\'{o}n, sobre todo si el autor no ha agregado una introducci\'{o}n que justifique racionalmente la no-obviedad a priori, y la profundidad de las relaciones halladas en sus enunciados.

Un criterio de no-trivialidad usado para un teorema nuevo en los arbitrajes de las comunicaciones cient\'{\i}ficas de matem\'{a}tica, termina a veces limit\'{a}ndose a  \lq\lq medir la longitud\rq\rq \  de su demostraci\'{o}n, verificar que no est\'{e} artificialmente alargada, y a \lq\lq contar las dificultades exitosamente sorteadas\rq\rq \  durante la misma.  Si estas mediciones est\'{a}n por abajo de ciertas cotas, usualmente el  teorema es tildado de \lq\lq trivial\rq\rq \  y por lo tanto \lq\lq irrelevante\rq\rq \ , y no es publicado, aunque su enunciado sea racionalmente sorprendente e inesperado.  Esta postura conduce, a veces, a que las relaciones profundas que puedan existir  entre conceptos matem\'{a}ticos de origen independiente, o de subdiciplinas matem\'{a}ticas diferentes, pueden permanecer difundidas muy parcialmente, solo como curiosidades o ejercicios, por fuera de la bibliograf\'{\i}a cient\'{\i}fica arbitrada y reconocida.

Concluimos que, por un lado el concepto de no-trivialidad o no-obviedad del enunciado de un teorema, y su percepci\'{o}n de profundidad por parte de los matem\'{a}ticos durante la actividad de investigaci\'{o}n,  est\'{a}n estrechamente ligados a criterios racionales informales no cuantitativos de relevancia y a propiedades intr\'{\i}nsecas de las relaciones matem\'{a}ticas contenidas en el enunciado, en forma parcialmente independiente de la econom\'{\i}a o abundancia de recursos formales utilizados en la demostraci\'{o}n.  Pero por otro parad\'{o}jicamente, para la comunicaci\'{o}n, aceptaci\'{o}n y publicaci\'{o}n de resultados matem\'{a}ticos nuevos, la brevedad o sencillez de recursos formales necesarios para la demostraci\'{o}n son frecuentemente consideradas se\~{n}ales de trivialidad u obviedad del enunciado.

\section{ Criterios de no-redundancia l\'{o}gico-formal y de redundancia \'{o}ptima informal.}

La redundancia l\'{o}gico-formal global de una subdisciplina o rama matem\'{a}tica es la propiedad de que al menos uno de sus  enunciados se puede deducir de los dem\'{a}s. En este sentido, todo resultado \lq\lq nuevo\rq\rq \  de una subdisciplina matem\'{a}tica ya existente, es redundante, excepto si para demostrarlo se requiri\'{o} agregar definiciones o presupuestos nuevos,  consistentes pero no deducibles de la teor\'{\i}a existente hasta el momento. Por este motivo la no-redundancia l\'{o}gico-formal, en su concepci\'{o}n global, no es un criterio aplicable usualmente en la matem\'{a}tica para juzgar la relevancia de sus enunciados.

Los criterios de no-redundancia l\'{o}gico-formales que efectivamente se tienen en cuenta en un juicio de valor de enunciados matem\'{a}ticos son dos, que establecen condiciones l\'{o}gicas diferentes e independientes entre s\'{\i}. Ambos son parciales en vez de globales.

El primero, que llamaremos no-redundancia mutua l\'{o}gico-formal, est\'{a} restringido a las proposiciones condicionantes (los supuestos, que deben ser mutuamente consistentes) cuando forman parte de una \'{u}nica definici\'{o}n, de un mismo conjunto de axiomas, o de la hip\'{o}tesis de un mismo teorema. Es el siguiente  criterio: el subconjunto de supuestos no es mutuamente redundante cuando ninguna de las proposiciones que lo forman puede deducirse de las dem\'{a}s.  Es m\'{a}s un criterio de bondad del enunciado respectivo que de su relevancia. Es decir, aunque el enunciado contenga en sus supuestos redundancias mutuas l\'{o}gico-formales no-obvias (no-triviales), puede a\'{u}n ser relevante, aunque se le considere mejorable, no \'{o}ptimo.

El segundo criterio de no-redundancia l\'{o}gico-formal de un teorema, que llamaremos no-redundancia relativa a la tesis, tambi\'{e}n est\'{a} restringida a las proposiciones condicionantes de la hip\'{o}tesis de un mismo teorema. Pero en vez de considerar la no-redundancia mutua, considera la necesidad de cada proposici\'{o}n supuesta en la hip\'{o}tesis,  para que se cumpla la tesis. Precisamente, el criterio es el siguiente: los supuestos en la hip\'{o}tesis de un teorema son redundantes desde el punto de vista l\'{o}gico-formal con respecto a la tesis del mismo teorema, si retirando por lo menos una de las proposiciones supuestas en la hip\'{o}tesis, de las dem\'{a}s se puede deducir la misma tesis del teorema.

Expongamos un ejemplo de redundancia l\'{o}gico-formal mutua de los supuestos de un enunciado.  En la definici\'{o}n de los difeomorfismos de Anosov se establece como condici\'{o}n que el comportamiento din\'{a}mico de los subfibrados invariantes sea uniformemente hiperb\'{o}lico. Un teorema no-trivial, demostrado por ejemplo en (Bonati et als., pp. 287-293), establece que si los subfibrados son invariantes y uniformemente hiperb\'{o}licos, entonces son continuos. Por lo tanto, esa propiedad de continuidad ser\'{\i}a redundante si se la agregara a la definici\'{o}n de difeormorfismo de Anosov, o a las hip\'{o}tesis de un teorema que la presuponga.

Sin embargo, y aunque la continuidad de los subfibrados es una propiedad redundante respecto a las otras condiciones que definen la clase de difeomorfismos de Anosov, frecuentemente !`se la agrega a la definici\'{o}n o a las hip\'{o}tesis de los teoremas! Se la incluye racionalmente, aunque sea redundante desde el punto de vista l\'{o}gico-formal. Por un lado esa redundancia no es obvia y la continuidad de los subfibrados es una propiedad no-trivial. Y por otro lado es una propiedad esencial y necesaria en las demostraciones de muchas caracter\'{\i}siticas relevantes (din\'{a}micas, topol\'{o}gicas y erg\'{o}dicas) de los difeomorfismos de Anosov. Concluimos que el criterio de no-redundancia l\'{o}gico-formal puede no aplicarse en algunos casos porque predomina el criterio de no-trivialidad racional informal. A\'{u}n as\'{\i}, en esos casos se acostumbra salvaguardar la forma agregando una nota que advierte de la redundancia l\'{o}gico-formal mutua de las condiciones supuestas.

Consideremos ahora un ejemplo de no-redundancia l\'{o}gico-formal de cada supuesto en la hip\'{o}tesis de un teorema con respecto a su tesis. El siguiente  ejemplo no corresponde a la investigaci\'{o}n matem\'{a}tica actual, sino  a la ense\~{n}anza de la matem\'{a}tica universitaria. Pero de todas formas es ilustrativo de los criterios de evaluaci\'{o}n de redundancias de los supuestos respecto a la tesis durante los procesos de investigaci\'{o}n matem\'{a}tica. El conocido teorema de Picard (ver por ejemplo Sotomayor, \cite[p. 13]{19}) establece que existe y es \'{u}nica la soluci\'{o}n con dato inicial de la siguiente ecuaci\'{o}n diferencial ordinaria:
\begin{equation}
 \frac{dx}{dt} = F(x), \end{equation}
donde la funci\'{o}n F es continua y Lipchitziana. La demostraci\'{o}n cl\'{a}sica de este teorema tiene dos partes no triviales: una, la prueba de existencia,  y otra, la prueba de unicidad, que aunque no-trivial resulta relativamente sencilla basada en la demostraci\'{o}n previa de la existencia. Nos interesa observar que la parte  m\'{a}s relevante de la demostraci\'{o}n, la de existencia de soluci\'{o}n, no requiere en realidad de la hip\'{o}tesis de Lipschitz de la funci\'{o}n F. M\'{a}s precisamente, el teorema de Peano (ver por ejemplo Sotomayor, \cite[ p. 16]{19}) establece que existe la soluci\'{o}n con dato inicial de la ecuaci\'{o}n diferencial (1) cuando la funci\'{o}n F del segundo miembro es continua, aunque no sea Lipschitziana. Por lo tanto, si la tesis a demostrar fuera solamente la existencia de soluci\'{o}n, la hip\'{o}tesis de Lipschitz ser\'{\i}a redundante respecto a la tesis.

Pero la tesis del teorema de Picard afirma tambi\'{e}n la unicidad de la soluci\'{o}n. Y aunque la demostraci\'{o}n de esa unicidad sea formalmente f\'{a}cil de obtener a partir de la demostraci\'{o}n de existencia, la hip\'{o}tesis de Lipschitz no es redundante. En efecto, es cl\'{a}sico el ejemplo de la ecuaci\'{o}n diferencial $dx/dt = F(x)$ donde $F$ es la funci\'{o}n ra\'{\i}z cuadrada. Para ella existe soluci\'{o}n con dato inicial $x(0)= 0$, pero no es \'{u}nica. Ese ejemplo muestra la no-redundancia de la hip\'{o}tesis de Lipschitz respecto a la tesis de unicidad del teorema de Picard, y por lo tanto la relevancia de esta condici\'{o}n, independientemente de que solo se la utilice en la parte m\'{a}s sencilla de la demostraci\'{o}n.

En la mayor\'{\i}a de los casos, cuando no es obvia, puede resultar dif\'{\i}cil verificar la no-redundancia l\'{o}gico-formal de los supuestos de un teorema, tanto la mutua en las proposiciones de la hip\'{o}tesis, como la de cada una de esas proposiciones relativa a la tesis.  Es frecuente la aparici\'{o}n de nuevos teoremas no-triviales en la matem\'{a}tica actual,  que establecen que un cierto supuesto de la hip\'{o}tesis de un teorema previo era redundante en relaci\'{o}n a la tesis.

En general, un juicio de valor racional de un teorema nuevo verifica,  adem\'{a}s de la consistencia y la no-trivialidad, la siguiente  pre-condici\'{o}n para la no-redundancia l\'{o}gico-formal. Ella  no asegura esta \'{u}ltima, pero es necesaria para esta. La pre-condici\'{o}n consiste  en que ninguna de las proposiciones supuestas en la hip\'{o}tesis de un nuevo teorema sea omitida, u obviamente omitible, a lo largo de todos los pasos deductivos de su demostraci\'{o}n. Nos referimos a la demostraci\'{o}n dada, no necesariamente a cualquier otra demostraci\'{o}n que pueda encontrarse en el futuro para el mismo teorema. Sin embargo, aunque necesaria, esta verificaci\'{o}n por s\'{\i} sola no asegura la no-redundancia (tambi\'{e}n llamada necesidad) de ese supuesto particular, para obtener la misma tesis del teorema mediante alguna otra demostraci\'{o}n a\'{u}n no conocida, a menos que se d\'{e} un contraejemplo. Un tal contraejemplo, que prueba la no-redundancia de cada proposici\'{o}n de la hip\'{o}tesis respecto a la tesis de un teorema, es un caso particular en que se cumplen todos, excepto uno, los supuestos en la hip\'{o}tesis, no se cumple ese uno, y no se cumple la tesis del teorema.  Nos referimos solamente a la necesidad de cada supuesto hipot\'{e}tico por separado, asumiendo verdaderos los dem\'{a}s y la tesis del teorema.  Por lo tanto, un criterio de relevancia l\'{o}gico-formal de un teorema nuevo considera positivamente, pero no obligatoriamente, los contraejemplos que puedan suministrarse adem\'{a}s del enunciado del teorema, que prueben que ninguno de los supuestos en la hip\'{o}tesis es redundante respecto a la tesis.

Otra cuesti\'{o}n relacionada con la no-redundancia de los supuestos relativa la tesis de un teorema, en general m\'{a}s dif\'{\i}cil de evaluar que la que definimos arriba, es la necesidad de todos los supuestos de la hip\'{o}tesis a la vez, para que valga la tesis. M\'{a}s precisamente, la no-redundancia de las hip\'{o}tesis de un teorema respecto a la tesis es total, cuando las proposiciones de la tesis implican por deducci\'{o}n la verdad de todas las proposiciones de la hip\'{o}tesis a la vez. Esta cuesti\'{o}n ya no se denomina en matem\'{a}tica propiedad de no-redundancia de la hip\'{o}tesis, sino una propiedad mucho m\'{a}s fuerte, llamada caracterizaci\'{o}n de la tesis. Es, desde el punto de vista l\'{o}gico-formal, equivalente a la verdad del teorema rec\'{\i}proco.

Frecuentemente, exceptuando los contraejemplos triviales, es muy dif\'{\i}cil que un teorema nuevo sea enunciado por primera vez acompa\~{n}ado de contraejemplos que prueban la necesidad (es decir la no-redundancia respecto a la tesis) de cada uno de los supuestos en su hip\'{o}tesis. Una de las fuentes de conjeturas que son planteadas por los matem\'{a}ticos investigadores, e incorporadas a sus agendas de investigaci\'{o}n, es precisamente la pregunta de si alguna de las afirmaciones asumidas en la hip\'{o}tesis de un teorema ya conocido, es necesaria, es decir no-redundante con respecto a la tesis del teorema  (y a las dem\'{a}s proposiciones de su hip\'{o}tesis).

Por ejemplo, Pesin en 1977, \cite{14} formul\'{o} por primera vez que para todos los difeomorfismos de Anosov de clase $C^2$, las medidas llamadas f\'{\i}sicas (i.e. medidas de probabilidad que describen la estad\'{\i}stica de \'{o}rbitas t\'{\i}picas seg\'{u}n Lebesgue), verifican una igualdad para la entrop\'{\i}a, que luego fue  llamada f\'{o}rmula de Pesin de la entrop\'{\i}a. (Ma\~{n}\'{e}, \cite[p. 95, f\'{o}rmula 2]{12}). \'{E}sta es una igualdad matem\'{a}tica entre la entrop\'{\i}a m\'{e}trica y el promedio espacial del logaritmo de la tasa de dilataci\'{o}n. La demostraci\'{o}n de la f\'{o}rmula de Pesin usa en forma esencial el supuesto hipot\'{e}tico de que el difeomorfismo es de clase $C^2$. Es indudablemente una hip\'{o}tesis no-redundante para que esa demostraci\'{o}n del teorema funcione.  Sin embargo, recientemente se descubri\'{o} que es redundante desde el punto de vista l\'{o}gico-formal, con respecto a la tesis. M\'{a}s precisamente, no es necesario que el difeomorfismo de Anosov sea de clase $C^2$ para  que toda medida f\'{\i}sica satisfaga la f\'{o}rmula de Pesin de la entrop\'{\i}a (ver por ejemplo Cerminara et als. \cite[pp. 737-761]{3}).

?`Pierde con este nuevo resultado su relevancia, respecto a la f\'{o}rmula de Pesin, la vieja hip\'{o}tesis de regularidad $C^2$? En nuestra opini\'{o}n, a favor de la cual argumentaremos a continuaci\'{o}n, y seg\'{u}n los argumentos racionales informales de buena parte de los matem\'{a}ticos que investigan en la  teor\'{\i}a erg\'{o}dica diferenciable, la respuesta es \lq\lq No, no la pierde en absoluto\rq\rq \ .  En efecto, la hip\'{o}tesis de regularidad $C^2$, si bien no es necesaria para deducir que las medidas f\'{\i}sicas satisfacen la f\'{o}rmula de Pesin, es necesaria para que el procedimiento anterior de demostraci\'{o}n de esta f\'{o}rmula, v\'{\i}a la propiedad de continuidad absoluta de las foliaciones invariantes, sea v\'{a}lida. Por un lado, esta propiedad de continuidad absoluta es considerada relevante en s\'{\i} misma por la comunidad matem\'{a}tica especializada en teor\'{\i}a erg\'{o}dica diferenciable,  seg\'{u}n pr\'{a}cticamente todos los criterios racionales informales que discutimos y analizamos a lo largo de este trabajo. Por otro lado, la demostraci\'{o}n anterior de la f\'{o}rmula de Pesin muestra, como paso intermedio, no solo que las medidas f\'{\i}sicas la satisfacen, sino que estas medidas son absolutamente continuas respecto a la foliaci\'{o}n inestable. Y este resultado, que no ser\'{\i}a cierto sin la hip\'{o}tesis de regularidad $C^2$ (como muestra el contraejemplo de Robinson y Young 1980, \cite{16}), tiene relevancia en s\'{\i} mismo: se considera no solo una herramienta de demostraci\'{o}n, sino una propiedad no-trivial relevante del sistema din\'{a}mico.

El ejemplo anterior muestra otra vez que los criterios l\'{o}gico-formales de no-redundancia pueden ser, en la pr\'{a}ctica de la investigaci\'{o}n matem\'{a}tica actual, inocuos respecto a un juicio de relevancia de enunciados, cediendo su lugar a criterios racionales informales de relevancia.  Es una de las razones por las que nos afiliamos a la integraci\'{o}n de los criterios l\'{o}gico-formales con los otros criterios racionales, pero informales, en los juicios de valor del quehacer de la investigaci\'{o}n matem\'{a}tica.

En la siguiente discusi\'{o}n distinguiremos dos sitios donde la redundancia racional informal, en vez de la no-redundancia l\'{o}gico-formal,  puede incidir positivamente en un juicio de relevancia de los resultados: la del proceso creativo o de descubrimiento de los enunciados matem\'{a}ticos, y la del proceso de comunicaci\'{o}n de \'{e}stos.

Es frecuente que los matem\'{a}ticos se ocupen de la comprensi\'{o}n profunda de los resultados matem\'{a}ticos y procuren comunicar el aspecto significativo de los contenidos,  adem\'{a}s de exponer la l\'{o}gica formal deductiva de las demostraciones.  Como muchas actividades humanas no exactas, la comunicaci\'{o}n de la matem\'{a}tica (y aunque la matem\'{a}tica sea una ciencia exacta, su comunicaci\'{o}n no lo es), requiere redundancia en la transmisi\'{o}n de informaci\'{o}n, para que los errores de expresi\'{o}n e interpretaci\'{o}n puedan ser corregidos, o por lo menos detectados, y minimizar sus efectos negativos en la eficiencia de la comunicaci\'{o}n. Aunque usando s\'{\i}mbolos y notaci\'{o}n muy espec\'{\i}ficos, el lenguaje matem\'{a}tico es un lenguaje humano. Definamos la entrop\'{\i}a del error como la tasa de expansi\'{o}n de la cantidad de informaci\'{o}n que es diferente entre la informaci\'{o}n transmitida y la recibida (estas diferencias son los \lq\lq errores\rq\rq \  en la transmisi\'{o}n). \lq\lq A pesar que raramente es mostrado en modelos diagram\'{a}ticos del proceso de comunicaci\'{o}n, la redundancia - la repetici\'{o}n de elementos dentro de un mensaje que previene el error en la comunicaci\'{o}n de informaci\'{o}n - es el mayor ant\'{\i}doto a la entrop\'{\i}a\rq\rq \  (Gordon, \cite[secci\'{o}n \lq\lq Entropy\rq\rq \ , p\'{a}rrafo 2]{7}). La mayor parte de los lenguajes escritos y hablados, por ejemplo, son aproximadamente 50% redundantes. M\'{a}s precisamente, en el lenguaje humano escrito o hablado, se podr\'{\i}a reducir a un 50% la cantidad de elementos, s\'{\i}mbolos o c\'{o}digos transmitidos, sin alterar la cantidad de informaci\'{o}n que se aspira comunicar. Aunque si se efectuara esa reducci\'{o}n, la seguridad por eventuales errores en la comunicaci\'{o}n ser\'{\i}a nula.

La redundancia natural en todo lenguaje humano (y quiz\'{a}s el lenguaje matem\'{a}tico sea de los que tienen menor redundancia entre todos los lenguajes humanos) es el mismo fen\'{o}meno f\'{\i}sico de la redundancia artificial de los medios de transmisi\'{o}n de la ingenier\'{\i}a moderna, que se dise\~{n}a para aumentar la seguridad de las comunicaciones digitales. Y tambi\'{e}n es el mismo fen\'{o}meno f\'{\i}sico de la redundancia exhibida en numerosos procesos naturales. Por ejemplo, en el an\'{a}lisis matem\'{a}tico estad\'{\i}stico de series temporales de datos, para estudiar fen\'{o}menos naturales como el estado del tiempo entre otras aplicaciones, se llama redundancia a la diferencia de la suma de las entrop\'{\i}as (definida la entrop\'{\i}a, en este caso, como el crecimiento de la cantidad de informaci\'{o}n significativa, y no como crecimiento de los errores de la informaci\'{o}n) de todas las variables por separado, menos la entrop\'{\i}a conjunta de esas variables. Dicho de otra forma, la cantidad de informaci\'{o}n de la serie de datos no es la suma de la informaci\'{o}n de cada dato por separado, sino que es esta suma menos la redundancia. Esta redundancia es mayor cuando m\'{a}s mutuamente dependiente son los datos estudiados. Por lo tanto, a igual cantidad de informaci\'{o}n de los datos por separado, cuanto mayor es la redundancia, m\'{a}s reducida es la cantidad de informaci\'{o}n significativa que se obtiene de ellos. Pero es justamente esta redundancia entre los datos, lo que permite hacer predicciones matem\'{a}ticas: cuanto mayor redundancia, m\'{a}s previsible ser\'{a} la evoluci\'{o}n futura de un fen\'{o}meno natural estudiado o artificial dise\~{n}ado.

Por razones de confiabilidad, se introduce naturalmente o artificialmente cierta cantidad de redundancia, la cual implica un sobre-dimensionamiento de los canales de transmisi\'{o}n o soporte de conservaci\'{o}n de esa informaci\'{o}n. Es decir, cuanto mayor es la redundancia, y por lo tanto mayor es la seguridad y fidelidad en la comunicaci\'{o}n y en la informaci\'{o}n transmitida,  menor es el aprovechamiento de los recursos necesarios para transmitirla. En efecto, por ejemplo en la teor\'{\i}a de las telecomunicaciones:

\em Mediante el proceso de conversi\'{o}n anal\'{o}gico-digital, cualquier medio disponible de telecomunicaciones tiene una capacidad limitada para la transmisi\'{o}n de datos. Esta capacidad es com\'{u}nmente medida por un par\'{a}metro llamado 'ancho de banda' . Como el ancho de banda de una se\~{n}al aumenta con el n\'{u}mero de bits a ser transmitidos, una funci\'{o}n importante en el sistema digital de comunicaciones es representar la se\~{n}al digitalizada con la menor cantidad posible de bits.  \rm (Lehnert,  \cite[ secci\'{o}n \lq\lq Source encoding\rq\rq \ , p\'{a}rrafo 1]{10}).

Pero para obtener eficiencia usando la menor cantidad posible de \lq\lq bits\rq\rq \  habr\'{\i}a que reducir a cero la redundancia en la cantidad de informaci\'{o}n de la se\~{n}al,  y por lo tanto la seguridad y confiabilidad de la transmisi\'{o}n de esta informaci\'{o}n. Esto muestra que en comunicaci\'{o}n de informaci\'{o}n, existe un compromiso entre la seguridad-confiabilidad en la transmisi\'{o}n de esa informaci\'{o}n, y la eficiencia en el uso del medio de transmisi\'{o}n utilizado. Al aumentar la redundancia se aumenta la primera pero se disminuye la segunda. Por lo tanto la situaci\'{o}n \'{o}ptima, seg\'{u}n sea el objetivo del que comunica, se encuentra en un valor intermedio de redundancia, que  no es nulo pero tampoco 100\%.

 En la comunicaci\'{o}n de resultados matem\'{a}ticos, aparece este mismo compromiso: as\'{\i} por ejemplo art\'{\i}culos matem\'{a}ticos que contienen resultados muy relevantes y son publicados en las revistas cient\'{\i}ficas m\'{a}s importantes, frecuentemente son m\'{a}s largos que el promedio de los art\'{\i}culos matem\'{a}ticos publicados en otras revistas cient\'{\i}ficas.  Esto sucede no solo porque la cantidad de enunciados de cada art\'{\i}culo excepcionalmente relevante sea relativamente grande, y sus demostraciones largas y complicadas, sino porque tambi\'{e}n contienen - usualmente en la introducci\'{o}n y en  notas remarcadas a lo largo del art\'{\i}culo - explicaciones sobre el significado profundo de esos enunciados y su relaci\'{o}n con otros resultados previamente conocidos.  Estas explicaciones, desde el punto de vista l\'{o}gico-formal son redundantes, pero a veces son imprescindibles en la b\'{u}squeda de efectividad y calidad de la comunicaci\'{o}n matem\'{a}tica.

Resulta parad\'{o}jico que, exceptuando casos como los descritos arriba, en buena parte de la bibliograf\'{\i}a de investigaci\'{o}n cient\'{\i}fica en matem\'{a}tica de hoy en d\'{\i}a, a pesar de que el recurso de transmisi\'{o}n de informaci\'{o}n utilizado (por ej. publicaci\'{o}n en l\'{\i}nea, en vez de papel) es abundante y de relativo bajo costo, algunos de los art\'{\i}culos sean escuetos, excesivamente concisos, y con poca redundancia en la comunicaci\'{o}n. Esto quiz\'{a}s se deba a que, para enunciados matem\'{a}ticos no excepcionalmente relevantes,  no sea tan necesaria la redundancia, pues cuando m\'{a}s f\'{a}cil de comprender o menos profundo es un enunciado, menos sujeto est\'{a} a errores en su transmisi\'{o}n o comunicaci\'{o}n.  La concisi\'{o}n, evitando la redundancia tanto formal-l\'{o}gica como informal en la comunicaci\'{o}n de resultados de investigaci\'{o}n, es tambi\'{e}n coherente con una postura filos\'{o}fica minimalista, frecuente entre algunos matem\'{a}ticos.

Previo a un juicio de no-redundancia de los supuestos de un enunciado, frecuentemente el investigador matem\'{a}tico procesa su comprensi\'{o}n profunda de esos supuestos a trav\'{e}s de un argumento por condicionantes contrafactuales. Estos son condicionantes que suponen la negaci\'{o}n de proposiciones asumidas previamente. En el ejemplo del teorema de Picard expuesto antes, si un estudiante tuviera que investigar la no-redundancia de la hip\'{o}tesis de Lipschitz, el condicionante contrafactual que asumir\'{\i}a es \lq\lq Si F no fuera Lipschitz ...\rq\rq \ . En ese caso la conclusi\'{o}n de no redundancia de la hip\'{o}tesis de Lipschitz establecer\'{\i}a: \lq\lq entonces la soluci\'{o}n con dato inicial podr\'{\i}a no ser \'{u}nica\rq\rq \ . Sin embargo esa conclusi\'{o}n no se obtiene por deducci\'{o}n l\'{o}gico-formal (en este ejemplo). Bajo el condicionante contrafactual,  es falso que la soluci\'{o}n con dato inicial sea necesariamente no \'{u}nica, como podr\'{\i}a esperar el estudiante .Es decir, no es en este caso solo un proceso deductivo el que lleva el condicionante contrafactual a la conclusi\'{o}n de no-redundancia de la hip\'{o}tesis que niega.

Concluimos que la discusi\'{o}n de los supuestos contrafactuales racionales, puede trascender la argumentaci\'{o}n l\'{o}gico-formal, y corresponde a una concepci\'{o}n filos\'{o}fica de racionalidad informal.
En el ejemplo que estamos discutiendo, la prueba de no-redundancia de la hip\'{o}tesis de Lipschitz en el teorema de Picard, requiere exhibir un contraejemplo para el cual se verifique ese condicionante, pero para el cual la soluci\'{o}n no sea \'{u}nica. Justamente, imaginar o descubrir ese contraejemplo, no es un resultado que se obtiene aplicando leyes l\'{o}gico-formales \'{u}nicamente. No se construye el contrajemplo meramente por deducci\'{o}n, sino que requiere de la creaci\'{o}n e imaginaci\'{o}n del investigador. \lq\lq La habilidad humana de pensar racionalmente sobre situaciones hipot\'{e}ticas y relaciones condicionantes se apoya en la capacidad de imaginar posibilidades\rq\rq \   (J. Laird, citado en Byrne \cite[p. 441]{2})

?`Pero es la creaci\'{o}n e imaginaci\'{o}n en matem\'{a}tica un proceso racional? Efectivamente lo es, afili\'{a}ndonos a la definici\'{o}n de imaginaci\'{o}n racional de Byrne, 2007, \cite{2}. \lq\lq En el pasado, racionalidad e imaginaci\'{o}n eran vistas como opuestos. Pero la investigaci\'{o}n ha mostrado que el pensamiento racional es m\'{a}s imaginativo de lo que se supon\'{\i}a\rq\rq \  (Byrne, \cite[p. 439]{2}. En particular, la imaginaci\'{o}n racional por condicionantes contrafactuales en matem\'{a}tica, permite crear alternativas en la b\'{u}squeda de nuevo conocimiento en la subdisciplina que se investiga. Aunque el proceso de imaginaci\'{o}n es informal, en la creaci\'{o}n de conocimiento matem\'{a}tico es racional. En efecto, \lq\lq los mismos principios del pensamiento racional, conducen tambi\'{e}n al pensamiento imaginativo... y el puente entre la racionalidad y la imaginaci\'{o}n puede ser construido sobre condicionales contrafactuales\rq\rq \  (Byrne \cite[p. 441]{2}. Sin embargo, aunque racional, la redundancia entre conjuntos de enunciados requerida durante los procesos de creaci\'{o}n e imaginaci\'{o}n en los que la investigaci\'{o}n matem\'{a}tica se apoya, no sigue un algoritmo l\'{o}gico-formal descriptible paso a paso. Por lo tanto, esa redundancia, y la manera en la que la imaginaci\'{o}n matem\'{a}tica se basa en ella,  son conceptos racionales informales.

Concluimos que la actividad de creaci\'{o}n e imaginaci\'{o}n racional del investigador matem\'{a}tico le permite formular conjeturas, ejemplos, contraejemplos, y tentativas de enunciados sobre los que construye y act\'{u}a su investigaci\'{o}n. En ese proceso, la redundancias l\'{o}gico-formales e informales constituyen un ingrediente esencial. No lo son solamente por las propiedades psicol\'{o}gicas de la creaci\'{o}n e imaginaci\'{o}n humanas, sino tambi\'{e}n por las caracter\'{\i}sticas intr\'{\i}nsecas del proceso de investigaci\'{o}n racional  de la matem\'{a}tica, y por las caracter\'{\i}sticas profundas de la matem\'{a}tica misma. M\'{a}s precisamente, del proceso de creaci\'{o}n e imaginaci\'{o}n racional en matem\'{a}tica, se aspira encontrar relaciones no triviales entre enunciados previos y/o nuevos. Es justamente la b\'{u}squeda de redundancias en esa informaci\'{o}n, que es anterior a los juicios de verdad en matem\'{a}tica, el camino que conduce la investigaci\'{o}n. Pero si el conjunto de enunciados sobre los que se investiga no fueran redundantes, esas relaciones buscadas no existir\'{\i}an, debido a la definici\'{o}n misma de independencia o no-redundancia de la cantidad de informaci\'{o}n de cada uno.  Si bien el investigador matem\'{a}tico usualmente no toma conciencia ni obtiene una medici\'{o}n precisa de la redundancia requerida entre los enunciados que considera en sus procesos de creaci\'{o}n e imaginaci\'{o}n racional, es en base a ella que puede construir el conocimiento matem\'{a}tico nuevo.

\section{Conclusi\'{o}n}

Expusimos y discutimos varias definiciones, ilustradas con ejemplos, de los conceptos de consistencia, no-trivialidad y redundancia en la investigaci\'{o}n de la matem\'{a}tica actual, tanto desde el punto de vista l\'{o}gico-formal como del racional informal.  Discutimos esos conceptos,  argumen-tando a favor de una postura, en cuanto a los criterios que inciden en los juicios de relevancia matem\'{a}tica, acorde con una filosof\'{\i}a de racionalidad formal e informal  integrada.

Aunque quiz\'{a}s no tanto como en \'{a}reas cient\'{\i}ficas basadas en la experimentaci\'{o}n y en las evidencias emp\'{\i}ricas, en la matem\'{a}tica tambi\'{e}n est\'{a} presente la racionalidad informal complementando a los criterios de racionalidad l\'{o}gico-formales, para integrar los juicios de relevancia de sus enunciados.

%%%%%%%%%%%%%%%%%%%%%%%%%%%%%%%%%%%%%%%%%%%%%%%%%%%%%%%%%%%%%%%%
%%%%%%%%%%%%%%%%%%%%%%%%%%%%%%%%%%%%%%%%%%%%%%%%%%%%%%%%%%%

\end{document}